# THE MOTION OF A SECOND CLASS PARTICLE FOR THE TASEP STARTING FROM A DECREASING SHOCK PROFILE

By Thomas Mountford and Hervé Guiol

*École Polytechnique Fédérale de Lausanne and Institut National Polytechnique de Grenoble*

We prove a strong law of large numbers for the location of the second class particle in a totally asymmetric exclusion process when the process is started initially from a decreasing shock. This completes a study initiated in Ferrari and Kipnis [*Ann. Inst. H. Poincaré Probab. Statist.* **13** (1995) 143–154].

**1. Introduction.** The totally asymmetric exclusion process (or TASEP) is an interacting system of indistinguishable particles on the one-dimensional lattice $\mathbb{Z}$. Each element of the lattice is called a *site*. Initially particles are distributed on $\mathbb{Z}$ according to the *exclusion rule* that prohibits multi-occupancy of sites (i.e., each site is at most occupied by one particle). Particles attempt to move one step (to their neighboring site) in a unique given direction (say to the right) independently at random times (exponential law with mean 1) provided this does not violate the exclusion rule. Two classic physical interpretations of the TASEP are in use: (i) as a moving interface on the plane (space–time); and (ii) as a toy model for traffic on a single-lane highway. The former gives a powerful tool to analyze the process in terms of a last passage percolation. This approach will be fully commented and developed throughout the paper. However, we postpone this to the following sections since we believe the second interpretation gives a more pleasant way to understand the result for the nonspecialist. We will thus begin to expose the results in this setting in an intuitive way.

1.1. *The traffic model.* We interpret particles as cars on a single-lane highway ($\mathbb{Z}$) with no possibility of passing (exclusion rule). For each car moving times are modeled by independent Poisson processes with mean 1. Due to the nonpassing rule, a car at position $x$ that attempts to move needs free space in front of it to occupy position $x+1$; otherwise the car just stays at $x$ beginning to jam, and so on....

Now suppose we look to this highway from far away, say from a helicopter (*macroscopic level*). As the mean velocity (jump rate of the process) is the same for all cars, the classical results about exclusion processes (see [9]) show that any given traffic density on the full highway is preserved over time.

Now typically if one starts with two different densities (which we call an initial *shock condition*), say $\lambda$ to the left of the origin and $\rho$ to the right, then two kinds of situations are possible. If $\lambda < \rho$, cars from the left are entering in a more jammed zone with density $\rho$. And what we see from our helicopter is that the jam front is moving to the left (which is called the *propagation of the shock*). When $\lambda > \rho$, then the jammed zone is to the left and what we observe in between is a resorption (i.e., *rarefaction fan* or rarefaction wave) of this jam into the (more) fluid region. Those facts are rigorously proved by the study of the *hydrodynamic* behavior of the TASEP (see, e.g., [12]).

Some natural questions then arise:

(a) In the first case $\lambda < \rho$: Is the shock sharp? That is, is there a microscopic identifier of the shock? The answer is yes. It has been proved by Ferrari, Kipnis and Saada [5] and then by Ferrari [3] that the *second class particle* identifies the shock (we refer to the previous articles for a rigorous statement of these facts). At the intuitive level the description of a second class particle is the following: It is a special particle that obeys all the rules described before except that other particles are insensitive to it; that is, if it stands just in front of a (regular) particle that tries to move, then it has to exchange its position with the moving particle. As one can see, a second class particle might move back, which is impossible for the other particles in this totally asymmetric context.

(b) In the second case $\lambda > \rho$: What is the microscopical counterpart of the resorption of the jam (rarefaction fan)? Ferrari and Kipnis [4] proved, first, that the position of the second class particle converges in distribution to a uniform law in the rarefaction fan, and second, that the second class particle, once having chosen an allowable velocity in the wedge (the rarefaction fan), remains close to it forever in probability. It was an open question to prove whether or not the above mentioned convergence holds more strongly. We prove indeed in the following that it is the case as a strong law of large numbers.

From now on we will abandon the traffic interpretation and switch to the interface model. To do so we start with the formal statement of the result and the mathematical description of the TASEP. We will then make the connection between the TASEP and the interface model and discuss some key ideas of the proof. We postpone the plan of the paper to the end of the next section.

### 1.2. *The interface settings.*

1.2.1. *The result.* We consider the position $X(t)$ at time $t \geq 0$ of a second class particle initially at the origin, that is, $X_0 = 0$, for a TASEP $(\eta_t)_{t \geq 0}$ for which at time 0, particles are independently present at sites $x, x \in \mathbb{Z} \setminus \{0\}$, with probability $\lambda$ for $x < 0$ and probability $\rho$ for $x > 0$.

We address (following [4]) the case $\lambda > \rho$, that is to say, an initial decreasing shock. The case $\lambda < \rho$ has already been completely resolved by the above mentioned papers [3] and [5]; see also [11] and [16] for similar results for more general processes.

The starting point is the result of Ferrari and Kipnis [4] that

$$(1) \qquad \text{as } t \to \infty \qquad \frac{X(t)}{t} \xrightarrow{D} U([1 - 2\lambda, 1 - 2\rho]),$$

where $U(I)$ denotes the uniform distribution on interval $I$; and

$$(2) \qquad \text{for } 0 < s < t \text{ fixed as } \varepsilon \to 0 \qquad \frac{X(t/\varepsilon)}{t/\varepsilon} - \frac{X(s/\varepsilon)}{s/\varepsilon} \xrightarrow{pr} 0.$$

Given this result it is natural to conjecture the following result, which will be proven in Section 4.

THEOREM 1. *For $(X(t))_{t \geq 0}$ as above there exists a uniform random variable $U$ on $[1 - 2\lambda, 1 - 2\rho]$ so that*

$$as \ t \to \infty \qquad \frac{X(t)}{t} \xrightarrow{a.s.} U.$$

For us the key ingredients are *Seppäläinen's variational formula* for TASEP (see [14] and (7) below) and *concentration inequalities* originating with Talagrand; see, for example, [1]. We do not need to use the exciting new results on last passage percolation of Johansson [6].

1.2.2. *The TASEP and its hydrodynamics.* The TASEP is an interacting particle system on $\{0, 1\}^{\mathbb{Z}}$ with generator on cylinder functions $f$

$$Lf(\eta) = \sum_{x \in \mathbb{Z}} \eta(x)[1 - \eta(x+1)](f(\eta^{x,x+1}) - f(\eta))$$

where

$$\eta^{x,x+1}(y) = \eta(y) \qquad \text{for } y \neq x \text{ or } x+1,$$
$$\eta^{x,x+1}(x) = \eta(x+1) \quad \text{and} \quad \eta^{x,x+1}(x+1) = \eta(x).$$

We interpret $\eta_t(x) = 1$ to mean that for the configuration $\eta_t$ there is a particle at site $x$ [or for the process $(\eta_s)_{s\geq 0}$ there is a particle at site $x$ at time $t$]. Particles try to move at exponential times to the site one to the right of their present site but moves to sites already occupied by another particle are suppressed. For details of more general exclusion processes see [9]; for details on Seppäläinen's description see [14, 15].

It is well known (see [10] and [12]) that the TASEP has the following scaling property. Let $(\eta_t^N)_{t\geq 0}$ for $N \in \mathbb{N}$ be a sequence of TASEPs such that for any finite interval $I \subset \mathbb{R}$

$$\left| \frac{1}{N} \sum_{x/N \in I} \eta_0^N(x) - \int_I u_0(r)\, dr \right| \xrightarrow{pr} 0,$$

where $u_0$ is a measurable function on $\mathbb{R}$ such that $0 \leq u_0(x) \leq 1$. Then for all $t \geq 0$

$$\left| \frac{1}{N} \sum_{x/N \in I} \eta_{Nt}^N(x) - \int_I u_t(r)\, dr \right| \xrightarrow{pr} 0,$$

where $u_t(r)$ is the *entropy solution* (see, e.g., [8]) of the scalar conservation law

$$(3) \qquad \frac{\partial u}{\partial t} + \frac{\partial G(u)}{\partial r} = 0$$

with flux function $G(u) = u(1-u)$ and initial condition $u_0$.

In particular, when $u_0(x) = \lambda \mathbf{1}_{x \leq 0} + \rho \mathbf{1}_{x > 0}$ with $\lambda > \rho$ the entropy solution produces a rarefaction fan

$$u_t(x) = \begin{cases} \lambda, & \text{if } x \leq (1-2\lambda)t, \\ (t-x)/2t, & \text{if } (1-2\lambda)t < x \leq (1-2\rho)t, \\ \rho, & \text{if } x > (1-2\rho)t. \end{cases}$$

Another way of having a look at this is to consider the (integrated) *Hamilton–Jacobi problem* (see, e.g., [2])

$$\frac{\partial U}{\partial t} + G\left(\frac{\partial U}{\partial r}\right) = 0$$

with $U_0$ satisfying, for all $x < y$,

$$U_0(y) - U_0(x) = \int_x^y u_0(r)\, dr.$$

Then the unique viscosity solution $U_t(x)$ of this problem is given by the *Hopf–Lax formula*

$$(4) \qquad U_t(x) = \sup_{y \in \mathbb{R}} \left\{ U_0(y) - t g\left(\frac{x-y}{t}\right) \right\},$$

where $g$ is the nonincreasing, nonnegative convex function such that for $u \in [0,1]$

$$G(u) = \inf_r \{ur + g(r)\},$$

that is, $g$ is the Legendre convex conjugate of the flux $G$. Note that the supremum in (4) is indeed achieved at some $y \in [x-t, x+t]$.

The solution $U_t$ is related to the entropy solution to the original equation (3) by the relation

$$\forall t \geq 0, \forall x < y \in \mathbb{R} \qquad U_t(y) - U_t(x) = \int_x^y u_t(w) \, dw.$$

Here and subsequently $g(x) = (1-x)^2/4$.

1.2.3. *Seppäläinen's variational formula.* Seppäläinen's formula gives a microscopic equivalent of (4) for the TASEP. We will now describe that formula; all the details can be found in [14].

We need first to introduce a *tool process* from which the TASEP can easily be retrieved. Let $(z_t)_{t \geq 0}$ be a server process on $\mathbb{Z}^{\mathbb{Z}}$, where $z_t(i)$ represents the position of the $i$th server of a system at time $t$. We impose the following exclusion rule:

$$0 \leq z_t(i+1) - z_t(i) \leq 1, \tag{5}$$

that is, two consecutive servers cannot overpass each other nor be too far (two sites or more apart).

The construction of the $z_t$ process is achieved by a system of independent Poisson processes. Let $\{(\mathcal{P}_i(t))_{t \geq 0}\}_{i \in \mathbb{Z}}$ be a collection of mutually independent Poisson processes with rate 1 on $]0, \infty[$, and call it a *Harris system*. At any epoch $\tau$ of $(\mathcal{P}_i(t))_{t \geq 0}$, $z_\tau(i)$ will be reduced by one unit provided this does not violate (5), in which case nothing happens to the system.

Given such a Harris system and an (independent) initial distribution $z_0 \in \mathbb{Z}^{\mathbb{Z}}$ that satisfies (5) on the same probability space we can construct the $z_t$ process at any time $t \geq 0$. The exclusion process is then retrieved via

$$\eta_t(x) = z_t(x) - z_t(x-1).$$

So the condition (5) is seen to be simply equivalent to the condition that $\eta_t(x) \in \{0, 1\}$.

Now we need to define a family $\{(w_t^k)_{t \geq 0} : k \in \mathbb{Z}\}$ of auxiliary processes, on the same probability space, such that each $(w_t^k)_{t \geq 0}$ is a server process like $(z_t)_{t \geq 0}$ satisfying the exclusion rule (5). Initially we define

$$w_0^k(i) = \begin{cases} z_0(k), & \text{if } i \geq 0, \\ z_0(k) + i, & \text{if } i < 0, \end{cases}$$

that is, all the servers with nonnegative label occupy the same position $z_0(k)$ and the others are put at distance 1 from their neighbors; dynamically

$w^k(i)$ attempts to jump to $w^k(i) - 1$ at the epochs of $(\mathcal{P}_{i+k}(t))_{t \geq 0}$.

The utility of the $w_t^k$ processes comes from the following *variational coupling* formula:

LEMMA 2 ([14], Lemma 4.1). *For all $i \in \mathbb{Z}$ and $t \geq 0$*,
$$z_t(i) = \sup_{k \in \mathbb{Z}} w_t^k(i - k) \qquad a.s. \tag{6}$$

The r.v.'s $\{w_t^k(i)\}_{\{-\infty < i < \infty\}}$ can be visualized as the height of an interface over the sites $i$. In order to start initially from level zero and obtain a growing surface, the family of interface processes $\{(\xi_t^k)_{t \geq 0} : k \in \mathbb{Z}\}$ is defined by

$$\xi_t^k(i) = z_0(k) - w_t^k(i) \qquad \text{for } i \in \mathbb{Z}, t \geq 0;$$

then

$$\xi_0^k(i) = \begin{cases} 0, & \text{if } i \geq 0, \\ -i, & \text{if } i < 0, \end{cases}$$

and the variational formula (6) gives Seppäläinen's variational formula

$$z_t(i) = \sup_{k \in \mathbb{Z}} \{z_0(k) - \xi_t^k(i - k)\}. \tag{7}$$

Observe that as the process $\xi^k$ does not depend on the initial $z_0(k)$, and depends on $k$ only through a translation of the indexing of $\{(\mathcal{P}_i(t))_{t \geq 0}\}$, dynamically

$\xi^k(i)$ jumps to $\xi^k(i) + 1$ at epochs of $\mathcal{P}_{i+k}$

provided the inequalities $\xi^k(i) \leq \xi^k(i - 1)$ and $\xi^k(i) \leq \xi^k(i + 1) + 1$ are not violated. Seppäläinen's variational formula is our key tool in the sense that it permits us to (see especially the proof of Proposition 19) trace the position of the second class particle.

1.2.4. *Strategy of the proof.* Loosely speaking, the strategy of the proof is the following: Keeping in mind Ferrari–Kipnis weak law (1), we want to take advantage of the idea contained in their other result (2), that is, once chosen a given velocity in the rarefaction fan the second class particle keeps following it.

To do so we analyze the trajectory of the second class particle by chopping it off into a sequence of increasing time intervals of order $(2^n)_{n \in \mathbb{N}}$. On each of these intervals we control the deviations from the original direction taken by the second class particle. This analysis is performed in Section 4 thanks to:

(a) Large deviation bounds obtained from a related last passage percolation problem described below (see Corollary 12 in Section 2), where the key tools are the above mentioned concentration inequalities and a trick from [7].

(b) A nice approximation of the server process by the solution of the Hamilton–Jacobi problem related with the hydrodynamic limit of the TASEP (Proposition 16 in Section 3). For this part the key tools are Hopf–Lax and Seppäläinen formulas.

Once obtained the almost sure convergence of $X(t)/t$ the weak law (1) suffices to conclude.

1.2.5. *The last passage percolation problem.* Following Seppäläinen [15] we now recall how formula (7) can be analyzed in terms of a last passage percolation problem.

We need first to introduce some notation. Define the *wedge* of admissible lattice paths denoted by

$$\mathcal{L} = \{(i,j) \in \mathbb{Z}^2 : j \geq 1, i \geq -j+1\},$$

with boundary $\partial \mathcal{L} = \{(i,0) : i \geq 0\} \cup \{(i,-i) : i < 0\}$.

For $(i,j) \in \mathcal{L} \cup \partial \mathcal{L}$ let

$$\Theta^k(i,j) = \inf\{t \geq 0 : \xi_t^k(i) \geq j\}$$

be the first time the interface $\xi^k$ reaches level $j$ at site $i$. The previous rules give

$$\Theta^k(i,j) = 0 \qquad \text{for } (i,j) \in \partial \mathcal{L}$$

and for $(i,j) \in \mathcal{L}$

$$\Theta^k(i,j) = \max\{\Theta^k(i-1,j), \Theta^k(i,j-1), \Theta^k(i+1,j-1)\} + \beta_{i,j}^k,$$

where $\beta_{i,j}^k$ is an exponential, mean 1, waiting time independent of the other $\beta_{i',j'}^k$.

Now consider the following last passage model: Let $\{t_{i,j} : (i,j) \in \mathcal{L}\}$ be a collection of i.i.d. exponential rate 1 random variables. Define the *passage times* $\{T(i,j)\}$ by

$$T(i,j) = 0 \qquad \text{for } (i,j) \in \partial \mathcal{L}$$

and

(8) $$T(i,j) = \max_{\pi \in \Pi(i,j)} \sum_{(m,\ell) \in \pi} t_{m,\ell} \qquad \text{for } (i,j) \in \mathcal{L}$$

where $\Pi(i,j)$ is the set of *admissible lattice paths*

$$\pi = \{(0,1) = (i_1,j_1), (i_2,j_2), \ldots, (i_p,j_p) = (i,j)\}$$

such that $(i_m, j_m) - (i_{m-1}, j_{m-1}) = (1,0)$ or $(-1,1)$.

Let
$$\xi_t(i) = \min\{j : (i, j+1) \in \mathcal{L},\ T(i, j+1) > t\},$$

with $\xi_0(j) = 0$ for $j \geq 0$ and $\xi_0(j) = -j$ for $j < 0$. Then from [14], *the process $\xi.(\cdot)$ has the same distribution as the process $\xi_.^k(\cdot)$ of* (7).

Furthermore, Seppäläinen [13] obtained

(9)
$$\frac{1}{n}\xi_{nt}([nx]) \xrightarrow{pr} tg(x/t) \equiv t\frac{(1-x/t)^2}{4} \qquad \text{for } -t \leq x \leq t,$$
$$\lim_{n\to\infty} \frac{1}{n} T([nx], [ny]) = \Gamma(x, y) \equiv (\sqrt{y} + \sqrt{x+y})^2;$$

here and in the sequel $[u]$ denotes the integer part of $u \in \mathbb{R}$. The limiting "shape" $g$ for $\xi.(\cdot)$ satisfies $\Gamma(x, g(x)) = 1$ for all $|x| \leq 1$, meaning it is a curve level of $\Gamma$.

1.3. *Organization of the paper.* The paper is organized as follows: in Section 2 we use a last passage percolation argument to obtain some simple, nonoptimal, large deviation bounds. In Section 3 these bounds are employed to show that for all $t \in [2^n/2, 2 \cdot 2^n]$, $\eta_t$ will be $2^{n\alpha}$ "close" to the hydrodynamic limit outside probability $\exp(-2^{n(1-\alpha)})$ for $\alpha$ close to 1 but strictly below it. In the final section an argument is given to show that a.s. $X(t)/t$ converges.

In the sequel we will use the following notation: $\mathbb{N} = \{1, 2, \ldots\}$ denotes the set of positive integers, $\mathbb{Z}^+ = \mathbb{N} \cup \{0\}$ denotes the set of nonnegative integers and $]u, v[$ (resp. $[u, v[$ or $]u, v])$ will denote the open interval (resp. semi-open intervals) with endpoints $u$ and $v$.

**2. Large deviation bounds.** In analyzing $\{T(i, j+1) < t\}$ or the a.s. equal event that $\xi_t(i) > j$ we consider the "longest" admissible path, $\pi$ [in the sense of the passage times of (8)] from $(0,0)$ to $(i, j+1)$ that passes through the lattice point $(0, 1)$ where $(i, j+1) \in \mathcal{L} = \{(k, \ell) \in \mathbb{Z}^2 : k + \ell \geq 1, \ell \geq 1\}$.

In this section we will consider a collection $\{\tau_{i,j} : (i, j) \in \mathbb{Z}^+ \times \mathbb{Z}^+\}$ of i.i.d., exponential mean 1, random variables and last passage time $\mathcal{T}(i, j)$ for $(i, j) \in \mathbb{Z}^+ \times \mathbb{Z}^+$ will be redefined as

(10)
$$\mathcal{T}(i, j) = \max_{\pi \in \mathbf{\Pi}(i,j)} \sum_{v=0}^{i+j-1} \tau_{\pi(v)} \qquad \text{for } (i, j) \in \mathbb{Z}^+ \times \mathbb{Z}^+,$$

where $\mathbf{\Pi}(i, j)$ is the set of *up-right admissible paths* from $(0, 0)$ to $(i, j)$ [starting at $(0,0)$]; that is, if $\pi \in \mathbf{\Pi}(i, j)$, then
$$\pi = (\pi(0) = (0, 0), \pi(1), \ldots, \pi(i+j) = (i, j))$$
where $\pi(v+1) - \pi(v) \in \{(1, 0), (0, 1)\}$ for $v \geq 1$.

To relate this to the previous section [and indeed the previous definition of $T(i,j)$], we are just making use of the isomorphism $(x,y) \in \mathcal{L} \mapsto (x+y-1, y-1) \in \mathbb{Z}^+ \times \mathbb{Z}^+$. It will be easy to obtain results for the original $T(i,j)$ and therefore the objects $\xi_t^y(x-y)$ from bounds on the redefined $\mathcal{T}(i,j)$'s. From (9) for $\theta \in (0,1)$ fixed,

$$\lim_{n \to \infty} \frac{\mathcal{T}([n\theta], [n(1-\theta)])}{n} = (\sqrt{\theta} + \sqrt{1-\theta})^2.$$

Let $\mathcal{V}_n := \{(i,j) \in \mathbb{Z}^+ \times \mathbb{Z}^+ : i+j = n\}$. The object of this section is to prove

PROPOSITION 3. *There exists $\varepsilon > 0$ so that for all $n$ large and $(i,j) \in \mathcal{V}_n$,*

$$\mathbb{P}(|\mathcal{T}(i,j) - n(\sqrt{\theta} + \sqrt{1-\theta})^2| \geq n^{1-\varepsilon}) \leq \exp(-n^\varepsilon),$$

*where $i = [n\theta]$ (and so $j = n - [n\theta]$).*

Our approach is to first obtain via concentration inequalities bounds on $\mathcal{T}(i,j) - \mathbb{E}\mathcal{T}(i,j)$ and then to consider $\mathbb{E}\mathcal{T}(i,j) - n(\sqrt{\theta} + \sqrt{1-\theta})^2$, for $n = i+j$ large. Our first result is:

PROPOSITION 4. *Let $(i,j) \in \mathcal{V}_n$. For any $x \geq 4$*

$$\mathbb{P}(|\mathcal{T}(i,j) - \mathbb{E}\mathcal{T}(i,j)| \geq 3\sqrt{n}\log(n^2)x)$$
$$\leq 2\left(2\exp\left(-\frac{(x-4)^2}{4}\right) + e\exp(-2\sqrt{n}\log(n^2)x)\right).$$

PROOF. First, following Kesten [7], consider the quantity

$$V_n = \sum_{(k,\ell) \in [0,n[^2 \cap \mathbb{Z}^2} (\tau_{k,\ell} - \log(n^2))^+.$$

Since $\tau_{k,\ell}$ are i.i.d. $\mathcal{E}xp(1)$ random variables we have for $t < 1$

$$\mathbb{E}\exp(tV_n) = \left(1 + \frac{1}{n^2}\frac{t}{1-t}\right)^{n^2} \leq \exp\left(\frac{t}{1-t}\right).$$

And so for any $y > 0$,

$$\mathbb{P}(V_n \geq y) \leq \inf_{0 < t < 1} \exp\left(\frac{t}{1-t} - ty\right) = e^{1-y}.$$

Now consider $\mathcal{T}'(i,j)$ derived from $\mathcal{T}(i,j)$ with $\tau_{\pi(i)}$ replaced by the bounded r.v. $\tau_{\pi(i)} \wedge 2\log n$.

Obviously $\mathcal{T}'(i,j) \leq \mathcal{T}(i,j) \leq \mathcal{T}'(i,j) + V_n$.

We analyze $\mathcal{T}'(i,j)$ using pages 62–64 of [1]. Let $M'_n$ be the median of $\mathcal{T}'(i,j)$.

Let $A = \{\omega : \mathcal{T}'(i,j)(\omega) \leq M'_n\}$. Suppose $\varpi$ is such that

$$\mathcal{T}'(i,j)(\varpi) \geq M'_n + \log(n^2)\sqrt{n}x;$$

then there exists an admissible up-right path $\pi$ from $(1,1)$ to $(i,j)$ [starting at $(0,0)$] such that

$$\sum_{i=0}^{n-1} \tau_{\pi(i)}(\varpi) \wedge \log(n^2) \geq M'_n + \log(n^2)\sqrt{n}x.$$

Then for any $\omega \in A$ we have

$$\sqrt{n}x \leq \sum_{i=0}^{n-1} \left|\left(\frac{\tau_{\pi(i)}(\varpi)}{\log(n^2)} \wedge 1\right) - \left(\frac{\tau_{\pi(i)}(\omega)}{\log(n^2)} \wedge 1\right)\right|$$

$$\leq \sum_{i=0}^{n-1} \mathbf{1}_{\tau_{\pi(i)}(\varpi) \neq \tau_{\pi(i)}(\omega)}.$$

Then with $\beta_k = 1/\sqrt{n}$, $k = 0,\ldots,n-1$, it follows that

$$g(A,\varpi) := \sup_{\{\beta:\|\beta\|_2 \leq 1\}} \inf_{\omega \in A} \sum_{k=0}^{n-1} \beta_k \mathbf{1}_{\tau_{\pi(i)}(\varpi) \neq \tau_{\pi(i)}(\omega)} \geq x,$$

where $\|\beta\|_2$ means the Euclidean norm of $\beta \in \mathbb{R}^n$. Thus applying Corollary 2.4.31 of [1] to $\mathbf{R}(A) = \mathbb{P}(\mathcal{T}'(i,j) \leq M'_n) \geq 1/2$, where $\mathbf{R}$ is the restriction of $\mathbb{P}$ on $\Sigma = [0, 2\log n]^{n^2}$, while

$$\mathbb{P}(\mathcal{T}'(i,j) \geq M'_n + \log(n^2)\sqrt{n}x) \leq \mathbf{R}(\{\varpi : g(A,\varpi) \geq x\})$$

we get

$$\mathbb{P}(\mathcal{T}'(i,j) \geq M'_n + \log(n^2)\sqrt{n}x) \leq 2\exp\left(-\frac{x^2}{4}\right).$$

Similarly

$$\mathbb{P}(\mathcal{T}'(i,j) \leq M'_n - \log(n^2)\sqrt{n}x) \leq 2\exp\left(-\frac{x^2}{4}\right).$$

Thus

$$\mathbb{E}\mathcal{T}'(i,j) \leq M'_n + \sqrt{n}\log(n^2) 2 \int_0^\infty \exp\left(-\frac{x^2}{4}\right) dx$$

$$\leq M'_n + 2\sqrt{\pi n}\log(n^2).$$

And similarly

$$\mathbb{E}\mathcal{T}'(i,j) \geq M'_n - 2\sqrt{\pi n}\log(n^2).$$

Thus we have

$$\mathbb{P}(\mathcal{T}'(i,j) \geq \mathbb{E}\mathcal{T}'(i,j) + \sqrt{n}\log(n^2)x)$$
$$\leq \mathbb{P}(\mathcal{T}'(i,j) \geq M'_n + \sqrt{n}\log(n^2)(x - 2\sqrt{\pi}))$$

and since $\sqrt{\pi} \leq 2$, for $x \geq 4$, $\leq 2\exp(-\frac{(x-4)^2}{4})$.

Thus since $\mathbb{E}\mathcal{T}'(i,j) \leq \mathbb{E}\mathcal{T}(i,j) \leq \mathbb{E}\mathcal{T}'(i,j) + e$, we have

$$\mathbb{P}(\mathcal{T}(i,j) \geq \mathbb{E}\mathcal{T}(i,j) + 3\sqrt{n}\log(n^2)x)$$
$$\leq \mathbb{P}(\mathcal{T}'(i,j) \geq \mathbb{E}\mathcal{T}(i,j) + \sqrt{n}\log(n^2)x) + \mathbb{P}(V_n \geq 2\sqrt{n}\log(n^2)x)$$
$$\leq 2\exp\left(-\frac{(x-4)^2}{4}\right) + \exp(1 - 2\sqrt{n}\log(n^2)x)$$

and get a similar bound for $\mathbb{P}(\mathcal{T}(i,j) \leq \mathbb{E}\mathcal{T}(i,j) - 3\sqrt{n}\log(n^2)x)$. $\square$

Next we concentrate on getting a useful bound for $|\mathbb{E}\mathcal{T}(i,j) - (\sqrt{i} + \sqrt{j})^2|$. We first assemble some elementary lemmas.

LEMMA 5. *Consider the random variable*

$$W = \sup_{(i,j) \in \mathcal{V}_n} \mathcal{T}(i,j).$$

*There exists a finite $K$ so that for any $c > K$ and all $n$ sufficiently large, for an event, $A$, of probability $\exp(-cn/(\log n)^2)$*

$$\mathbb{P}(W \geq cn|A) \leq \exp\left(-\frac{cn}{4}\right).$$

PROOF. For an up-right path of length $n$ starting at $(0,0)$ there are $2^n$ possible choices. For such a path, say $\pi$, the probability that

$$\sum_{i=0}^{n-1} \tau_{\pi(i)} \geq cn$$

is $\leq \exp(-cn/2)V^n$ where $V := \mathbb{E}\exp(X/2)$ and $X \sim \mathcal{E}xp(1)$ so

$$\mathbb{P}(W \geq cn|A) \leq \exp\left(\frac{cn}{(\log n)^2}\right) 2^n V^n \exp\left(-\frac{cn}{2}\right) \leq \exp\left(-\frac{cn}{4}\right)$$

for $c$ sufficiently large. $\square$

LEMMA 6. *For all positive $x, y$ and positive integer $n$,*

$$\mathbb{E}\mathcal{T}([nx], [ny]) \leq n(\sqrt{x} + \sqrt{y})^2.$$

PROOF. Suppose $\mathbb{E}\mathcal{T}([nx],[ny]) > n(\sqrt{x}+\sqrt{y})^2$. For each integer $L$
$$[Lnx] \geq L[nx] \quad \text{and} \quad [Lny] \geq L[ny];$$
and so $\mathcal{T}([Lnx],[Lny]) \geq \mathcal{T}(L[nx], L[ny])$.

But the longest path from $(0,0)$ to $(L[nx], L[ny])$ is longer than the longest path from $(0,0)$ to $(L[nx], L[ny])$ which goes through points $(\ell[nx], \ell[ny])$ for $0 \leq \ell \leq L$. This second quantity is equal to
$$\sum_{\ell=1}^{L} Z_\ell$$
where $Z_\ell$, are i.i.d. r.v.'s with
$$\mathbb{E}Z_\ell = \mathbb{E}\mathcal{T}([nx],[ny]).$$
Thus by the strong law of large numbers we have a.s.
$$\liminf_{L\to\infty} \frac{\mathcal{T}([Lnx],[Lny])}{nL} \geq \mathbb{E}Z_1 > (\sqrt{x}+\sqrt{y})^2.$$
But we have from [13]
$$\frac{\mathcal{T}([Lnx],[Lny])}{nL} \to (\sqrt{x}+\sqrt{y})^2$$
when $L \to \infty$ and this gives a contradiction and the lemma follows. $\square$

We also record a simple result that will be needed later.

LEMMA 7. *For all $\theta \in ]0, 1/2[$, $\theta > \delta > 0$, if $x_i \in ]\theta - \delta, \theta + \delta[^c$ for $i = 1, 2, \ldots, n$ and for $\sum \alpha_i = 1$, $\alpha_i \geq 0$, $\sum \alpha_i x_i = \theta$, then*
$$\sum \alpha_i (\sqrt{x_i} + \sqrt{1-x_i})^2 \leq (\sqrt{\theta}+\sqrt{1-\theta}) - 2\delta^2.$$

PROOF. For $f(x) = (\sqrt{x}+\sqrt{1-x})^2$
$$f''(x) = -\frac{1}{\sqrt{x}\sqrt{1-x}} - \frac{1}{2}\sqrt{\frac{1-x}{x^3}} - \frac{1}{2}\sqrt{\frac{x}{(1-x)^3}} \leq -4,$$
so
$$f(x) \leq f(\theta) + f'(\theta)(x-\theta) - 2(x-\theta)^2. \qquad \square$$

PROPOSITION 8. *If $\varepsilon > 0$ is sufficiently small, then for all $n$ sufficiently large it is the case that for each $\theta \in [0,1]$ there exist $(x,y) \in \mathcal{V}_n$ such that $\|(x,y) - (n\theta, n(1-\theta))\| \leq n^{4/5}$ and*
$$\mathbb{E}\mathcal{T}(x,y) \geq n(\sqrt{\theta}+\sqrt{1-\theta})^2 - n^{1-\varepsilon}.$$

PROOF. First if $\theta \leq n^{-1/10}$ (or by symmetry $\theta \geq 1 - n^{-1/10}$), then we can consider any deterministic path $\pi : (0,0) \to ([n\theta], n - [n\theta])$. Thus we can take, for example, $x = [n\theta]$. We have

$$\mathbb{E}\mathcal{T}([n\theta], n - [n\theta]) \geq \mathbb{E}\sum_{i=0}^{n-1} \tau_{\pi(i)} = n$$
$$= n(\sqrt{\theta} + \sqrt{1-\theta})^2 - 2n\sqrt{\theta}\sqrt{1-\theta}$$
$$\geq n(\sqrt{\theta} + \sqrt{1-\theta})^2 - n^{1-\varepsilon}$$

for $\varepsilon < 1/20$ and $n$ large enough. Hence it is enough to consider $\theta \in ]n^{-1/10}, 1 - n^{-1/10}[$. We fix $\varepsilon$ to be small. We assume there exists a $\theta$ and an $n$ for which the condition fails and obtain a contradiction if $n$ is sufficiently large (and $\varepsilon$ has been fixed to be sufficiently small).

Fix a relevant "direction" $\theta$. We suppose that $\theta \geq 1/2$ without loss of generality and that for all $(x,y) \in \mathcal{V}_n$ such that $\|(x,y) - (n\theta, n(1-\theta))\| \leq n^{4/5}$, we have

$$\text{(11)} \qquad \mathbb{E}\mathcal{T}(x,y) < n(\sqrt{\theta} + \sqrt{1-\theta})^2 - n^{1-\varepsilon}.$$

Now by Proposition 4 we have if (11) holds, then for all $(x,y) \in \mathcal{V}_n$

$$\text{(12)} \qquad \begin{aligned} &\mathbb{P}(\mathcal{T}(x,y) \geq n(\sqrt{\theta} + \sqrt{1-\theta})^2 - n^{1-\varepsilon}/2) \\ &\leq 2\exp\left(-\frac{n^{2(1-\varepsilon)}}{288n(\log n)^2}\right) := \frac{1}{M(n)} \equiv \frac{1}{M} \end{aligned}$$

and, also by Proposition 4 and Lemma 6, for all $(x,y) \in \mathcal{V}_n$ [not necessarily "close" to $(n\theta, n(1-\theta))$] we have

$$\text{(13)} \qquad \begin{aligned} \mathbb{P}(\mathcal{T}(x,y) \geq (\sqrt{x} + \sqrt{y})^2 + n^K) &\leq 2\exp\left(-\frac{n^{2K}}{288n(\log n)^2}\right) \\ &:= \frac{1}{N(n)} \equiv \frac{1}{N} \end{aligned}$$

for $K > 1/2$.

A suitable value for $K$ will be chosen later but for the moment we assume that $1/2 < K < 1 - \varepsilon$ and that $n$ is sufficiently large to ensure that $M \gg N$.

For each $(m,r) \in \mathbb{Z} \times \mathbb{Z}^+$ we say the block $B(m,r) := \{(u,v) \in \mathcal{V}_{rn} : u \in [nr/2 + n(m-1/2), nr/2 + n(m+1/2)[\}$ is *bad* (otherwise we say it is *good*) if there exists $(u,v) \in B(m,r)$ such that (at least) one large deviation event of type (12) or (13) holds for the system translated by $(u,v)$.

The probability that a block is good is at least

$$1 - 2\frac{n^2}{N}.$$

We are now in a position to sketch our approach to the proof.

14 T. MOUNTFORD AND H. GUIOL

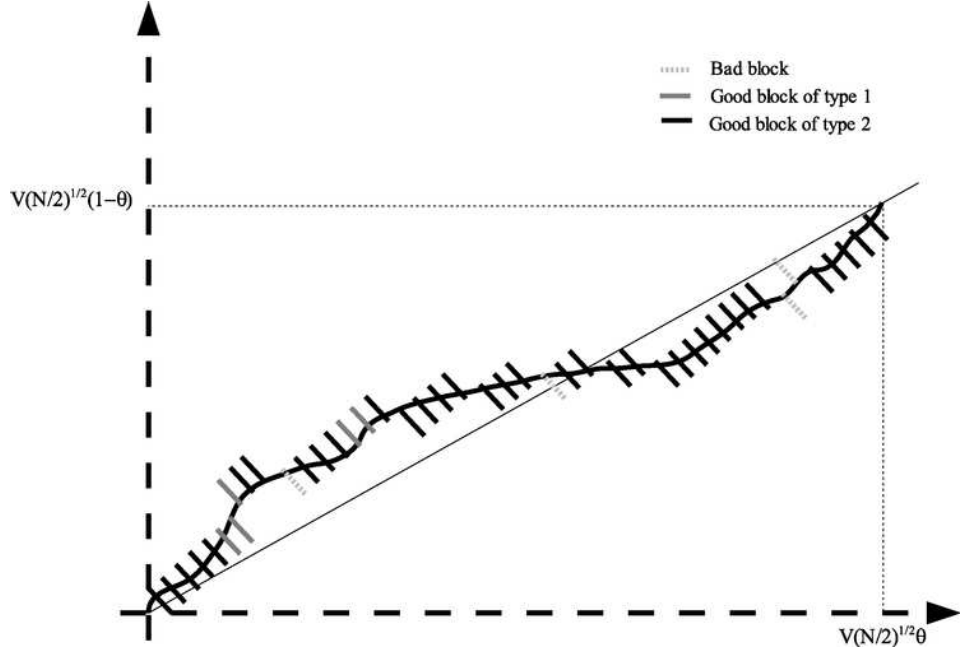

FIG. 1. *Representation of a path $\pi$ from $(0,0)$ to $V\sqrt{(N/2)}(\theta, 1-\theta)$.*

For any path $\pi$ from $(0,0)$ to $V\sqrt{(N/2)}(\theta, 1-\theta)$, we will consider it in subsegments of length $n$. Without loss of generality we may assume that $\sqrt{(N/2)}/n$ is an integer. We follow [7] and use the concentration inequalities to show that over all paths $\pi$ the contribution to the "length" of $\pi$ from $n$-segments of the path with great deviations is very small. This is summarized in Lemma 10 below. This will leave us two sorts of $n$-length path segments: those (type 2: see Figure 1) whose increment is "within" $n^{4/5}$ of $n(\theta, 1-\theta)$ and those (type 1) which are not. For the first collection and our assumption on the relevant expectations the average contribution should be less than $n(\sqrt{\theta} + \sqrt{1-\theta})^2$. For the latter segments we will only have Lemma 6 as a bound on the expectations but Lemma 7 will enable us again to conclude that the average contribution will fall short. The limit (9) will then be invoked to give a contradiction.

For calculation purposes let us by randomization have r.v.'s

$$\psi_{i,j} = \begin{cases} 1, & \text{with probability } 2n^2/N, \\ 0, & \text{with probability } 1 - 2n^2/N, \end{cases}$$

so that $\{B(i,j) \text{ is not good}\} \subseteq \{\psi_{i,j} = 1\}$ and for $(i_1, j_1), (i_2, j_2), \ldots, (i_r, j_r)$, $i_1 < i_2 < \cdots < i_r$, $\{\psi_{i_k, j_k}\}_{1 \leq k \leq r}$ are independent.

Now we renormalize by considering hyperblocks $G(i,j)$ [see Figure 2 where we have chosen $N = 8n^2$ for illustration purpose; note that in the sequel we

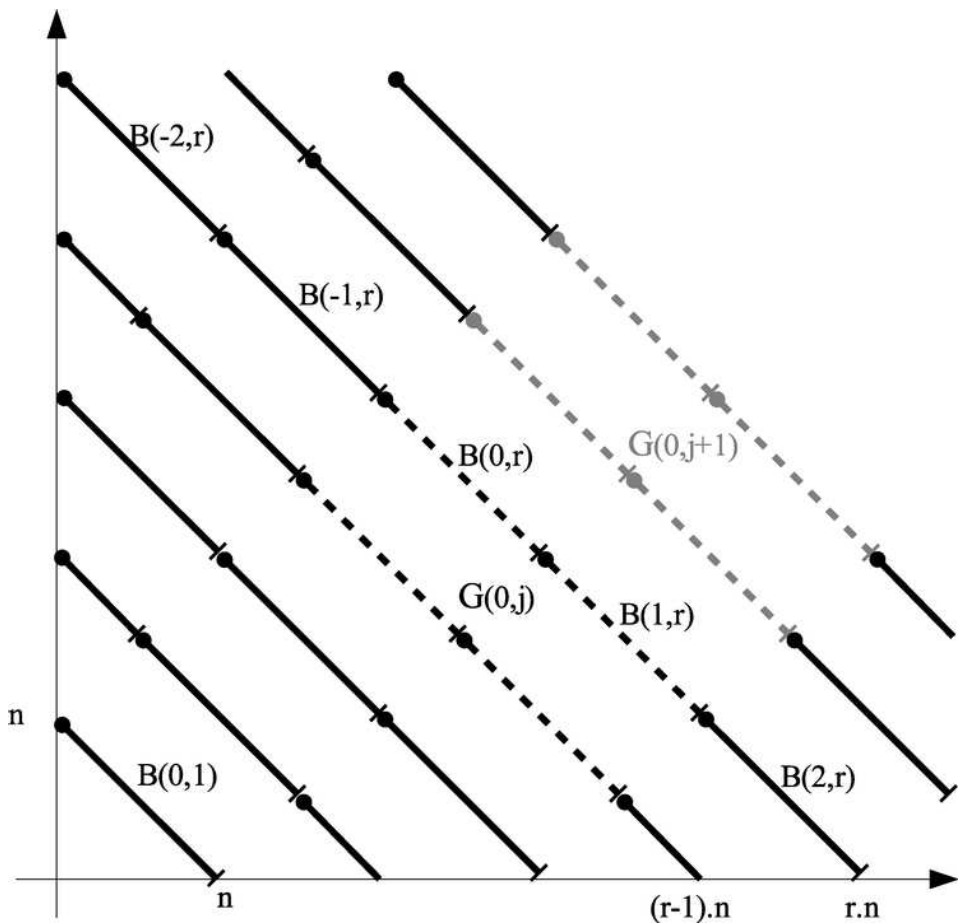

FIG. 2. *Blocks $B(\ell,r)$ and hyperblocks $G(i,j)$ (here $N = 8n^2$).*

will take $N$ of order $\exp(n^\alpha)$, with $\alpha > 0$] which are the union of

$$\bigcup_{\substack{i\sqrt{(N/2)}/n \leq \ell < (i+1)\sqrt{(N/2)}/n \\ j\sqrt{(N/2)}/n \leq r < (j+1)\sqrt{(N/2)}/n}} B(\ell, r).$$

The r.v. $X^G(i,j)$ is

$$\sum_{\substack{i\sqrt{(N/2)}/n \leq \ell < (i+1)\sqrt{(N/2)}/n \\ j\sqrt{(N/2)}/n \leq r < (j+1)\sqrt{(N/2)}/n}} W_{\ell,r}\psi_{\ell,r}$$

where
$$W_{\ell,r} = \sup\left\{\sum_{i=0}^{n-1}\tau_{\pi(i)}: \text{ over paths } \pi \text{ starting in } B(\ell,r)\right\}.$$

Then we have by independence and Lemma 5 that:

LEMMA 9. *For all $i,j$*
$$\mathbb{E}\exp\left(\frac{X^G(i,j)}{9n}\right) \leq H$$

*for constant $H$ not depending on $i,j$ or $n$.*

PROOF. We write $X^G(i,j)$ as the sum $X_o^G(i,j) + X_e^G(i,j)$ where
$$X_o^G(i,j) = \sum_{\substack{i\sqrt{(N/2)}/n \leq \ell < (i+1)\sqrt{(N/2)}/n \\ j\sqrt{(N/2)}/n \leq r < (j+1)\sqrt{(N/2)}/n \\ \ell \text{ odd}}} W_{\ell,r}\psi_{\ell,r}.$$

The reason for the introduction of these supplementary random variables is that all the terms in the sum defining $X_o^G(i,j)$ or $X_e^G(i,j)$ are independent, which is not the case for $X^G(i,j)$. By Cauchy–Schwarz one has
$$\mathbb{E}\exp\left(\frac{X^G(i,j)}{9n}\right) = \mathbb{E}\exp\left(\frac{X_o^G(i,j) + X_e^G(i,j)}{9n}\right)$$
$$\leq \left(\mathbb{E}\exp\left(\frac{2X_o^G(i,j)}{9n}\right)\right)^{1/2}\left(\mathbb{E}\exp\left(\frac{2X_e^G(i,j)}{9n}\right)\right)^{1/2}.$$

The lemma now follows from the claimed independence and Lemma 5. □

Thus for any path $\pi$ from $(0,0)$ to $v_1 + v_2 = V\sqrt{(N/2)}(\theta, 1-\theta)$ there is a corresponding *$G$-level path* $\pi^G: (0,0) \to (1,j_1) \to (2,j_2) \to \cdots \to (V,j_V)$. There are $2^V$ such paths and for all path $\pi^G$ we have an r.v.
$$Z(\pi^G) = \sum_{i=0}^{V-1} X^G(i,j_i).$$

By independence and the previous result
$$\mathbb{E}\exp\left(\frac{Z(\pi^G)}{9n}\right) \leq H^V$$

so
$$\mathbb{P}(Z(\pi^G) \geq n^3 V) \leq H^V \exp\left(-\frac{Vn^2}{9}\right),$$

and so

$$\mathbb{P}\left(\sup_{\pi^G} Z(\pi^G) \geq n^3 V\right) \leq H^V 2^V \exp\left(-\frac{Vn^2}{9}\right)$$

which is $\leq 2^{-V}$ if $V \geq V_0$.

We can in a similar (and easier) way redo this analysis with random variable $X^G(i,j)$ replaced by

$$\sum_{\substack{i\sqrt{(N/2)}/n \leq \ell < (i+1)\sqrt{(N/2)}/n \\ j\sqrt{(N/2)}/n \leq r < (j+1)\sqrt{(N/2)}/n}} \psi_{\ell,r}.$$

We have shown more than the following.

LEMMA 10. *With probability tending to 1, as $V \to \infty$ for all paths $\pi$ from $(0,0)$ to $V\sqrt{(N/2)}(\theta, 1-\theta)$*

$$\sum_{i=0}^{V\sqrt{(N/2)}/n-1} I_{\{\pi(in) \text{ is in a bad block}\}} \sum_{k=in}^{(i+1)n-1} \tau_{\pi(k)} \leq Vn^3$$

*and the number of bad blocks covered occurring in such a path is similarly bounded.*

PROOF. Let us denote by $A = A(V)$ the event that for all paths from $(0,0)$ to $V\sqrt{(N/2)}(\theta, 1-\theta)$ the above mentioned bounds hold.

Now fix a path $\pi$ and let vectors $e(i) = (e(i)_1, e(i)_2) := \pi(in) - \pi((i-1)n)$, $i = 1, 2, \ldots, V\sqrt{(N/2)}/n$. Let $V_1 =$ collection of $i$ so that $\pi((i-1)n)$ is in a good block and $\|e(i) - (n\theta, n(1-\theta))\| > n^{4/5}$. Let $V_2 =$ collection of $i$ so that $\pi((i-1)n)$ is in a good block and $\|e(i) - (n\theta, n(1-\theta))\| \leq n^{4/5}$.

By definition on the event $A$ for every path $\pi$,

$$\sum_{(r,s) \in \pi} \tau_{r,s} \, I_{\{(r,s) \text{ is in a bad block}\}} \leq Vn^3,$$

thus

$$(14) \quad \sum_{i=0}^{V\sqrt{(N/2)}-1} \tau_{\pi(i)} \leq Vn^3 + \sum_{i \in V_1} \sum_{k=(i-1)n+1}^{in} \tau_{\pi(k)} + \sum_{i \in V_2} \sum_{k=(i-1)n+1}^{in} \tau_{\pi(k)}.$$

The sum over $V_2$ is bounded above by

$$(15) \quad |V_2|\left(n(\sqrt{\theta} + \sqrt{1-\theta})^2 - \frac{n^{1-\varepsilon}}{2}\right) = |V_2|\left(nf(\theta) - \frac{n^{1-\varepsilon}}{2}\right)$$

by definition of $V_2$ and *goodness*, where $f(x) = (\sqrt{x} + \sqrt{1-x})^2$, while, again using the definition of $V_1$ and *goodness*, the sum over $V_1$ is bounded above by

$$\sum_{i \in V_1} (nf(v_i) + n^K) \tag{16}$$

where $v_i = n^{-1} e(i)_1$.

Since the sum of $e(i)$ over all $i$ must equal $V\sqrt{(N/2)}(\theta, 1-\theta)$ we have, from the second part of Lemma 10, that the sum of the terms $n^{-1} e(i)$ over $V_1$ and $V_2$ must be within $Vn^3$ of $n^{-1} V \sqrt{(N/2)}(\theta, 1-\theta)$ on event $A$. Furthermore by definition, the average of the sum of $n^{-1} e(i)$ over $V_2$ must be within $n^{-1/5}$ of $(\theta, (1-\theta))$. From these two facts we have that the average of the $v_i$ for $i \in V_1$ must equal $\theta + r$ where, first,

$$|r| \leq \frac{n^4}{\sqrt{(N/2)} - n^4} + \frac{|V_2|}{|V_1|}\left(\frac{1}{n^{1/5}} + \frac{\sqrt{(N/2)}}{\sqrt{(N/2)} - n^4}\right)$$

and, second, for $n$ large enough and $|V_1| \geq |V_2|$, it must be the case that

$$|r| \leq \frac{1}{2n^{1/10}}. \tag{17}$$

We also note that under our hypothesis for $n$ large and $i \in V_2$, provided $\varepsilon$ was picked sufficiently small

$$nf(v_i) \geq nf(\theta) - \frac{n^{1-\varepsilon}}{6}. \tag{18}$$

From (14), (15), (16) and (18) the sum over variables associated to path $\pi$ satisfies

$$\begin{aligned}
\sum_{i=0}^{V\sqrt{(N/2)}-1} \tau_{\pi(i)} &\leq Vn^3 + \sum_{i \in V_1} nf(v_i) + |V_1| n^K \\
&\quad + \sum_{i \in V_2} nf(v_i) - |V_2| \frac{n^{1-\varepsilon}}{3} \\
&= Vn^3 + \sum_{i \in V_1 \cup V_2} nf(v_i) + |V_1| n^K - |V_2| \frac{n^{1-\varepsilon}}{3}.
\end{aligned} \tag{19}$$

But under event $A(V)$ the average, $(|V_1| + |V_2|)^{-1} \sum_{i \in V_1 \cup V_2} nv_i$, is within $n^3/\sqrt{(N/2)}$ of $\theta$ and so by concavity of $f$ and elementary bounds on function $f$,

$$\sum_{i \in V_1 \cup V_2} nf(v_i) \leq (|V_1| + |V_2|)\left(nf(\theta) + L\frac{n^{3/2}}{(N/2)^{1/4}}\right),$$

where $L$ does not depend on $n, N$ and $\theta$. Therefore using (19) the sum over path $\pi$ satisfies

$$
\begin{aligned}
(20) \quad \sum_{i=0}^{V\sqrt{(N/2)}-1} \tau_{\pi(i)} &\leq Vn^3 + V\sqrt{(N/2)}f(\theta) \\
&\quad + Ln^{1/2}(N/2)^{1/4} + |V_1|n^K - |V_2|\frac{n^{1-\varepsilon}}{3} \\
&\leq V\sqrt{(N/2)}f(\theta) - V\sqrt{(N/2)}\frac{n^K}{n}
\end{aligned}
$$

provided $|V_2|/|V_1| \geq 4n^K/n^{1-\varepsilon}$, $1/2 < K < 1-\varepsilon$ and $n$ is sufficiently large uniformly on $\pi$.

It remains to deal with $\pi$ so that $|V_2|/|V_1| < 4n^K/n^{1-\varepsilon}$. In this case the average of the $v_i$ over $i \in V_1$ must be within $n^{4/5}/(2n)$ of $\theta$ and we have by Lemma 7 (with $\delta = n^{-1/5}/2$) and using (17) and the inequality $f(\theta+r) \leq f(\theta) + Cn^{1/10}|r|$, where $C$ does not depend on $N, n$ or $\theta$, that

$$
\sum_{i \in V_1} nf(v_i) \leq -|V_1|\frac{n^{2(4/5-1)+1}}{4} + C|V_1|\frac{|V_2|}{|V_1|}n^{9/10} + n|V_1|f(\theta)
$$

$$
= -|V_1|\frac{n^{3/5}}{4} + C|V_2|n^{9/10} + n|V_1|f(\theta).
$$

From this and, once more (14), (15) and (16), we have

$$
\sum_{i=0}^{V\sqrt{(N/2)}-1} \tau_{\pi(i)} \leq Vn^3 + V\sqrt{(N/2)}f(\theta) \\
+ |V_1|\left(\frac{-n^{3/5}}{4} + n^K\right) + |V_2|\left(\frac{-n^{1-\varepsilon}}{2} + Cn^{9/10}\right)
$$

which is $\leq V\sqrt{(N/2)}f(\theta) - V\sqrt{(N/2)}n^{-\varepsilon}/5$ if $3/5 > K > 1/2$ and $1-\varepsilon > 9/10$ and $n$ is sufficiently large. In either of these two cases we have, on event $A$ (whose probability tends to 1 as $V$ tends to infinity), for all admissible paths $\pi$

$$
\frac{1}{V\sqrt{(N/2)}} \sum_{i=0}^{V\sqrt{(N/2)}-1} \tau_{\pi(i)} \leq f(\theta) - \frac{n^{-\varepsilon}}{5}.
$$

That is, on event $A$

$$
\frac{\mathcal{T}([V\sqrt{(N/2)}\theta], [V\sqrt{(N/2)}(1-\theta)])}{V\sqrt{(N/2)}} \leq (\sqrt{\theta} + \sqrt{1-\theta})^2 - \frac{n^{-\varepsilon}}{5}
$$

and so $\mathcal{T}([V\sqrt{(N/2)}\theta],[V\sqrt{(N/2)}(1-\theta)])/(V\sqrt{(N/2)})$ does not tend in probability to $(\sqrt{\theta}+\sqrt{1-\theta})^2$. This contradiction establishes the desired result.

□

The proof of Proposition 8 is now complete.　□

The next result extends the preceding technical result to a more useful one.

COROLLARY 11. *There exists $0<\varepsilon_0<1/8$ so that for $n$ sufficiently large and for all $x,y \in \mathbb{N}$ such that $x+y=n$,*

$$\mathbb{E}\mathcal{T}(x,y) \geq (\sqrt{x}+\sqrt{y})^2 - n^{1-\varepsilon_0}.$$

PROOF. Fix $\varepsilon_1$ so small that Proposition 8 applies for all $m \geq n_0$ with $\varepsilon = \varepsilon_1$. Furthermore suppose that $\sqrt{n} \geq n_0 + 2$.

There are two cases to consider: $x \wedge y \leq n^{1-\varepsilon_1}$ and $x \wedge y > n^{1-\varepsilon_1}$.

We start with the second hypothesis: $x \wedge y > n^{1-\varepsilon_1}$ (since $x+y=n$, this implies that $\varepsilon_1 > 2/\log n$, but this will not be a problem for large enough $n$); without loss of generality we suppose that $\varepsilon_1 < 1/50$. Take $x_1 = [x/\sqrt{n}]$ and $y_1 = [y/\sqrt{n}]$; then obviously $\sqrt{n} - 2 \leq x_1 + y_1 \leq \sqrt{n}$. As

$$\theta := \frac{x_1}{x_1+y_1} > \frac{n^{1-\varepsilon_1}}{n} - \frac{1}{\sqrt{n}} = n^{-\varepsilon_1} - n^{-1/2},$$

thus $> n^{-1/10}$ in the chosen range for $\varepsilon_1$ provided $n_0$ was fixed sufficiently large. By Proposition 8 there exist $(x_2,y_2) \in \mathbb{N} \times \mathbb{N}$ with $x_2+y_2 = x_1+y_1$ and $\|(x_2,y_2)-(x_1,y_1)\| \leq (x_1+y_1)^{4/5} \leq n^{2/5}$ such that

$$\mathbb{E}\mathcal{T}(x_2,y_2) \geq (x_1+y_1)\left(\sqrt{\frac{x_1}{x_1+y_1}}+\sqrt{\frac{y_1}{x_1+y_1}}\right)^2 - (x_1+y_1)^{1-\varepsilon_1}$$

$$\geq (\sqrt{x_1}+\sqrt{y_1})^2 - n^{(1-\varepsilon_1)/2}.$$

Now consider the set of paths from $(0,0)$ to $(x,y)$ that pass through

$$(x_2,y_2),(2x_2,2y_2),\ldots,([\sqrt{n}(1-n^{-(\varepsilon_1-1/5)/2})]x_2,$$
$$[\sqrt{n}(1-n^{-(\varepsilon_1-1/5)/2})]y_2) \leq (x,y).$$

The expectation of the maximum of such paths is at least

$$\sqrt{n}(1-n^{-(1/10-\varepsilon_1/2)})\mathbb{E}\mathcal{T}(x_2,y_2)$$
$$\geq \sqrt{n}(1-n^{-(1/10-\varepsilon_1/2)})(\sqrt{x_1}+\sqrt{y_1})^2 - n^{(1-\varepsilon_1)/2}$$
$$\geq (\sqrt{x}+\sqrt{y})^2 - n^{1-\varepsilon_0}$$

if $\varepsilon_0$ was fixed sufficiently small.

We now treat the second case $x \wedge y \leq n^{1-\varepsilon_1}$. We suppose without loss of generality that $x \leq n^{1-\varepsilon_1}$. In this case $(\sqrt{x} + \sqrt{y})^2 \leq (\sqrt{n^{1-\varepsilon_1}} + \sqrt{y})^2 \leq y + 2n^{1-\varepsilon_1}$. We now consider the path from $(0,0)$ to $(0,y)$ and reason as in the start of Lemma 10. $\square$

This result together with Lemma 6 through Proposition 4 gives Proposition 3.

We finish the section by making the link to Seppäläinen's representation explicit.

COROLLARY 12. *For $\varepsilon < \varepsilon_0$, given in Corollary 11, and $(x,y) \in \mathcal{L} = \{(i,j) \in \mathbb{Z}^2 : j \geq 1, i \geq -j+1\}$ with $t = (\sqrt{x+y} + \sqrt{y})^2$ (i.e., in our notation $y = tg(x/t)$ and $x \in [-t,t]$),*

$$\mathbb{P}(|\xi_{nt}([nx]) - [ny]| > [tn]^{1-\varepsilon}) \leq 3\exp(-(nt)^\varepsilon).$$

PROOF. Without loss of generality we suppose that $(nt)^{1-\varepsilon}$ is an integer. The event $\{\xi_{nt}([nx]) > [ny] + (nt)^{1-\varepsilon}\}$ is equal to the event $\{\mathcal{T}([nx] + [ny] + (nt)^{1-\varepsilon} - 1, [ny] + (nt)^{1-\varepsilon} - 1) \leq nt\}$. Therefore

$$\mathbb{P}(\xi_{nt}([nx]) > [ny] + n^{1-\varepsilon})$$
$$= \mathbb{P}(\mathcal{T}([nx] + [ny] + (nt)^{1-\varepsilon} - 1, [ny] + (nt)^{1-\varepsilon} - 1) \leq nt).$$

Observe that the longest path from $(0,0)$ to $(z+k, w+k)$, $k, w, z \in \mathbb{N}$, would be bigger than the sum of the longest path from $(0,0)$ to $(z,w)$ and the passage time of the given path $\gamma = ((z,w), (z+1,w), \ldots, (z+k,w), (z+k, w+1), \ldots, (z+k, w+k))$. Thus

$$\mathbb{P}(\xi_{nt}([nx]) > [ny] + n^{1-\varepsilon})$$
$$\leq \mathbb{P}(\mathcal{T}([nx] + [ny], [ny]) \leq nt - (nt)^{1-\varepsilon})$$
$$+ \mathbb{P}\left(\sum_{i=1}^{(nt)^{1-\varepsilon}-1} \tau_{[nx]+[ny]+i,[ny]} \right.$$
$$\left. + \sum_{j=1}^{(nt)^{1-\varepsilon}-1} \tau_{[nx]+[ny]+(nt)^{1-\varepsilon}+1,[ny]+j} \leq (nt)^{1-\varepsilon}\right).$$

By Proposition 3, for $n$ sufficiently large the first probability is bounded by $\exp(-(tn)^\varepsilon)$, while by elementary large deviation bounds the second is bounded by

$$K\exp(-c(tn)^{1-\varepsilon}) < \frac{\exp(-(tn)^\varepsilon)}{2}$$

for $n$ large.

Arguing similarly

$$\begin{aligned}
&\mathbb{P}(\xi_{nt}([nx]) < [ny] - (nt)^{1-\varepsilon}) \\
&\leq \mathbb{P}(\mathcal{T}([nx] + [ny] - (nt)^{1-\varepsilon} - 1, [ny] - (nt)^{1-\varepsilon} - 1) \geq nt) \\
&\leq \mathbb{P}(\mathcal{T}([nx] + [ny], [ny]) \geq nt + (nt)^{1-\varepsilon}) \\
&\quad + \mathbb{P}\Bigg(\sum_{i=1}^{(nt)^{1-\varepsilon}-1} \tau_{[nx]+[ny]-(nt)^{1-\varepsilon}+i,[ny]-(nt)^{1-\varepsilon}} \\
&\quad + \sum_{j=1}^{(nt)^{1-\varepsilon}-1} \tau_{[nx]+[ny],[ny]-(nt)^{1-\varepsilon}+j} \leq (nt)^{1-\varepsilon}\Bigg).
\end{aligned}$$

And we obtain similar bounds for the probability of a large value of $\xi_{nt}([nx])$. □

REMARK. Given $t$, the above result immediately gives a probabilistic bound on $\xi_{nt}([nx])$ for $x$ where there exists $y$ such that $t = (\sqrt{x+y} + \sqrt{y})^2$, that is, for $x \in [-t, t]$. The usual bounds on Poisson processes allow us to deal with deviations of $\xi_{nt}([nx])$ for other $x$. For instance, if $x < -t$ one has

$$\mathbb{P}(|\xi_{nt}([nx]) + [nx]| \geq n^{1-\varepsilon}) \leq \mathbb{P}(\xi_{nt}([-nt - n^{1-\varepsilon}]) \geq 1 + [-nt - n^{1-\varepsilon}])$$

which can be bounded by the probability that a Poisson random variable of parameter $nt$ exceeds $[nt + n^{1-\varepsilon}]$.

**3. Hydrodynamic consequences at particle level.** For a configuration $\eta^n \in \{0,1\}^{\mathbb{Z}}$ indexed by integer $n$ and a measurable function $u_0$ taking values in $[0,1]$, we say

$$\eta^n \overset{M,v}{\sim} u_0$$

if for every $x \in [-Mn, Mn]$,

$$\left|\sum_{y=-Mn}^{x} \eta^n(y) - n\int_{M}^{x/n} u_0(r)\,dr\right| \leq v.$$

The main result of this section is:

PROPOSITION 13. *Let $M \geq 8$, $t \in \,]1/4, 4[$ and $\varepsilon < \varepsilon_0$, the constant of Corollary 11. For $(\eta_t)_{t\geq 0}$ an exclusion process with $\eta_0^n \overset{M,n^{1-\varepsilon}}{\sim} u_0$ and $n$ sufficiently large outside probability $\exp(-(nt)^{\varepsilon/2})$,*

$$\eta_{nt}^n \overset{M/2, 6n^{1-\varepsilon}}{\sim} u_t$$

*where $(u_s)_{s\geq 0}$ is the unique entropic solution to the scalar conservation law* (3) *with initial data $u_0$.*

Before beginning the proof of this result we give a simple lemma which enables us to reduce the analysis of general exclusion processes to that of finite systems.

LEMMA 14. *Let $(\eta_t)_{t\geq 0}$ be an exclusion process and for $k \in \mathbb{Z}_+$ let $(\eta_t^k)_{t\geq 0}$ be the exclusion process generated by the same Harris system as $(\eta_t)_{t\geq 0}$ and satisfying $\eta_0^k(x) = \eta_0(x) I_{|x|\leq k}$, $\forall x \in \mathbb{Z}$. For $M > 2, \exists c = c(M) > 0$ so that*

$$\mathbb{P}(\eta_s^{Mn}(x) = \eta_s(x) \ \forall |x| \leq Mn/2, 0 \leq s \leq n) \geq 1 - e^{-cn}.$$

PROOF. Let $(\eta_t^L)_{t\geq 0}$ [resp. $(\eta_t^R)_{t\geq 0}$] be exclusion processes (resp. reversed dynamics exclusion processes, i.e., $1 - p = 1$: total asymmetry to the left) run by the same Harris system as $(\eta_t)_{t\geq 0}$ so that a point $t \in \mathbb{P}_x$ represents a potential jump from site $x + 1$ to $x$ for process $(\eta_t^R)_{t\geq 0}$ and with initial configurations given by $\eta_0^L(x) = \delta_{x,-([Mn]+1)}$ and $\eta_0^R(x) = \delta_{x,[Mn]+1}$, where as before $[\cdot]$ denotes the integer part and $\delta$ is the Kronecker delta function. Then the event

$$\{\eta_s^{Mn}(x) = \eta_s(x) \ \forall |x| \leq Mn/2, 0 \leq s \leq n\}^c$$

is contained in the event

$$\left\{\sum_{x \geq -Mn/2} \eta_n^L(x) > 0\right\} \cup \left\{\sum_{x \leq Mn/2} \eta_n^R(x) > 0\right\}.$$

But the probability of these latter two events is simply equal to that of a Poisson random variable of parameter $n$ exceeding $[Mn] + 1 - [Mn/2]$ and the lemma follows. □

REMARK. This corresponds to a property of entropy solutions.

PROOF OF PROPOSITION 13. We first assume that $u_0 \equiv 0$ outside $[-M, M]$ and that $\eta_0^n(x) = 0$ for $x \notin [-Mn, Mn]$. Our condition gives that for $z_0^n(x) = \sum_{y=-Mn}^{x} \eta_0^n(y) = \sum_{y=-\infty}^{x} \eta_0^n(y)$ and $U_0(x) = \int_{-M}^{x} u_0(r)\,dr = \int_{-\infty}^{x} u_0(r)\,dr$,

(21) $$\forall x \quad \left| z_0^n(x) - nU_0\left(\frac{x}{n}\right) \right| \leq n^{1-\varepsilon}.$$

We have by the Hopf–Lax formula that $U_t(s) = \int_{-\infty}^{s} u_t(r)\,dr$ satisfies

$$U_t(s) = \sup_{r \in \mathbb{R}} \left\{ U_0(r) - tg\left(\frac{s-r}{t}\right) \right\},$$

and by Seppäläinen's formula for $x \in \mathbb{Z}$,

$$z_t(x) = \sup_{y \in \mathbb{Z}} \{z_0(y) - \xi_{x-y}^y(t)\}.$$

Let us for the moment fix $x \in [-Mn/2, Mn/2]$. From the observation after (4) we can choose $r^* \in [-(M+8)/2, (M+8)/2]$ so that

$$U_t\left(\frac{x}{n}\right) = \sup_{r \in \mathbb{R}}\left\{U_0(r) - tg\left(\frac{x-nr}{nt}\right)\right\} = U_0(r^*) - tg\left(\frac{x-nr^*}{nt}\right).$$

Take $y$ to be the integer part of $nr^*$. It is immediate from Seppäläinen's formula and (21) that

$$z_{nt}(x) \geq z_0(y) - \xi_{nt}^y(x-y) \geq nU_0\left(\frac{y}{n}\right) - n^{1-\varepsilon} - \xi_{nt}^y(x-y).$$

By Corollary 12 outside of probability $3\exp-(nt)^\varepsilon$ for $n$ large, we have that

$$\xi_{nt}^y(x-y) \leq ntg\left(\frac{x-y}{nt}\right) + n^{1-\varepsilon}$$

and so

$$z_{nt}(x) \geq nU_0\left(\frac{y}{n}\right) - 2n^{1-\varepsilon} - ntg\left(\frac{x-y}{nt}\right)$$

$$\geq nU_0(r^*) - 3n^{1-\varepsilon} - ntg\left(\frac{x-nr^*}{nt}\right)$$

$$= nU_t\left(\frac{x}{n}\right) - 3n^{1-\varepsilon}$$

by the Lipschitz properties of $U_0$ and $g$ and our choice of $r^*$. Thus for all $x \in [-Mn/2, Mn/2]$

(22) $$\mathbb{P}\left(z_{nt}(x) < nU_t\left(\frac{x}{n}\right) - 3n^{1-\varepsilon}\right) \leq 3\exp-(nt)^\varepsilon.$$

The argument for the converse is similar: By the finiteness assumption on $\eta_0^n$

$$z_{nt}(x) = \sup_{y \in \mathbb{Z}}\{z_0(y) - \xi_{nt}^y(x-y)\} = \sup_{|y| \leq Mn}\{z_0(y) - \xi_{nt}^y(x-y)\}.$$

So the event $\{z_{nt}(x) > nU_t(x/n) + 3n^{1-\varepsilon}\}$ is the union $\bigcup_{|y| \leq Mn}\{z_0(y) - \xi_{nt}^y(x-y) > nU_t(x/n) + 3n^{1-\varepsilon}\}$ and consequently

$$\mathbb{P}\left(z_{nt}(x) > nU_t\left(\frac{x}{n}\right) + 3n^{1-\varepsilon}\right)$$

$$\leq \sum_{|y| \leq Mn} \mathbb{P}\left(z_0(y) - \xi_{nt}^y(x-y) > nU_t\left(\frac{x}{n}\right) + 3n^{1-\varepsilon}\right).$$

We fix integer $y$ in the relevant range:

$$\mathbb{P}\left(z_0(y) - \xi_{nt}^y(x-y) > nU_t\left(\frac{x}{n}\right) + 3n^{1-\varepsilon}\right)$$

$$\leq \mathbb{P}\left(\xi_{x-y}^y(nt) < nU_0\left(\frac{y}{n}\right) - nU_t\left(\frac{x}{n}\right) - 2n^{1-\varepsilon}\right)$$

$$\leq \mathbb{P}\left(\xi_{nt}^y(x-y) < ntg\left(\frac{x-y}{nt}\right) - 2n^{1-\varepsilon}\right),$$

by the Hopf–Lax identity. But, having recourse once more to Corollary 12 and the remark following it, we bound this latter probability by $3\exp-(nt)^\varepsilon$ if $n$ is large.

After summing over $y$ we find that

$$\mathbb{P}\left(z_{nt}(x) > nU_t\left(\frac{x}{n}\right) + 3n^{1-\varepsilon}\right) \leq (3\exp-(nt)^\varepsilon)(2nM+1).$$

Summing this and (22) over $x$ in the relevant range, we have that outside probability (bounded by) $\exp-(nt)^{\varepsilon/2}$, if $n$ is sufficiently large,

$$\forall\, |x| \leq \frac{Mn}{2} \qquad \left|z_{nt}(x) - nU_t\left(\frac{x}{n}\right)\right| \leq 3n^{1-\varepsilon}.$$

The proof of the proposition is completed by appealing to Lemma 14 for the general case. $\square$

In the following lemma and subsequent proposition let $U_t(x)$ be the (integrated) solution to the Hamilton–Jacobi problem starting from the initial condition

$$U_0(x) = \begin{cases} \rho x, & \text{for } x > 0, \\ \lambda x, & \text{for } x < 0, \end{cases}$$

and let $(z_t(x))_{x\in\mathbb{Z}}$ be the server process associated to $(\eta_t)_{t\geq 0}$, the exclusion process beginning as product measure $\rho x$ for $x \geq 0$; product measure $\lambda x$ for $x < 0$ with $z_0(0) = 0$.

LEMMA 15. *For $\varepsilon_0 < 1/8$, $z_0$ and $U_0$ as above outside probability $\exp-2^{n/4}$ (for $n$ large)*

$$\left|z_0(x) - 2^n U_0\left(\frac{x}{2^n}\right)\right| \leq \frac{2^{n(1-\varepsilon_0)}}{10},$$

*for all $|x| \leq 8 \cdot 2^n$.*

PROOF. We fix $x$ positive and $\leq 8 \cdot 2^n$ without loss of generality. Denoting $\text{Bin}(n,p)$ a binomial random variable with parameters $n$ and $p$

$$z_0(x) \stackrel{D}{=} \text{Bin}(x,\rho) + 1 \qquad (\text{resp. } 0)$$

26    T. MOUNTFORD AND H. GUIOL...26    T. MOUNTFORD AND H. GUIOL

depending on whether we are considering first class particles only or first and second class together. In either case

$$P\left(\left|z_0(x) - U_0\left(\frac{x}{2^n}\right)2^n\right| > \frac{2^{n(1-\varepsilon_0)}}{10}\right) \leq \mathbb{P}\left(|\operatorname{Bin}(x,\rho) - x\rho| \geq \frac{2^{n(1-\varepsilon_0)}}{10} - 1\right)$$

$$\leq \exp(-2^{n(1-3\varepsilon_0)})$$

for $n$ large by standard bounds on binomial random variables (and using the bound $\varepsilon_0 < 1/8$). Similar bounds hold for $x$ negative and the result follows by summing. □

This lemma and Proposition 13 immediately give:

PROPOSITION 16.  *Let $(z_t(i))_{i \in \mathbb{Z}, t > 0}$ be the server process associated to the TASEP $(\eta_t)_{t \geq 0}$ defined in the Introduction. Let $U_t(x)$ be the (integrated) solution to the Hamilton–Jacobi problem starting from the initial condition $U_0(x) = \rho x$ for $x > 0$; and $\lambda x$ for $x < 0$. We suppose that*

$$z_0(0) = 0, \qquad U_0(0) = 0;$$

*then for $\varepsilon_1 < \varepsilon_0$ as defined in Corollary 11 or Proposition 13 and $\forall t \in {]2^n/2, 2 \cdot 2^n]}$ and each $y \in [-4 \cdot 2^n, 4 \cdot 2^n]$*

$$\left|z_t(y) - tU_t\left(\frac{y}{t}\right)\right| \leq t^{1-\varepsilon_1}$$

*outside probability $k \exp -t^{\varepsilon_1/2}$, provided $n$ is sufficiently large.*

**4. The a.s. convergence.**  In this section we wish to assemble the established results to prove Theorem 1.

Given Ferrari and Kipnis [4] it suffices to show that $X(t)/t$ converges a.s., from (1) the distribution of the limit random variable will follow immediately. This will be accomplished if we can show that for $\delta > 0$ arbitrarily small

$$\limsup_{t \to \infty} \frac{X(t)}{t} - \liminf_{t \to \infty} \frac{X(t)}{t} \leq \delta$$

with probability at least $1 - \delta$. We fix $\delta > 0$ and take integer $m$ (a power of 2) so that $1/m < \delta/10$. For $n$ positive and $i = 0, 1, 2, \ldots m$, we use $t_i^n$ to denote the time $2^n(1 + i/m)$ (so for all $n, t_m^n = t_0^{n+1}$). The following elementary result will enable us to restrict attention to $X(t)$ for $t$ equal to $t_i^n$ for some $n, i$.

LEMMA 17. *For $(X(t))_{t\geq 0}$ and $t_i^n$ as previously defined, a.s.*
$$\limsup_{n\to\infty} \sup_{0\leq i\leq m-1} \sup_{t_i^n\leq t\leq t_{i+1}^n} \left|\frac{X(t)}{t} - \frac{X(t_i^n)}{t_i^n}\right| \leq \frac{2}{m} < \frac{\delta}{5}.$$

PROOF. For the second class particle we associate two rate-1 Poisson processes: $N_+$ which contains $t$ if and only if $t \in \mathcal{P}_{X(t-)}$ (i.e., $t$ corresponds to a forward jump time of the second class particle) and $N_-$ which contains $t$ if and only if $t \in \mathcal{P}_{X(t-)-1}$ (i.e., $t$ corresponds to a backward jump time of the second class particle). The event $\{\sup_{t_i^n\leq t\leq t_{i+1}^n} X(t) - X(t_i^n) > (1+\varepsilon)2^n/m\}$ is contained in the event $\{N_+(t_{i+1}^n) - N_+(t_i^n) > (1+\varepsilon)2^n/m\}$ which has probability bounded by $K\exp(-c2^n\varepsilon)$ by elementary large deviations for Poisson random variables. A similar situation holds for the event $\{\sup_{t_i^n\leq t\leq t_{i+1}^n} X(t_i^n) - X(t) > (1+\varepsilon)2^n/m\}$. Thus we have by the Borel–Cantelli lemma that for each $\varepsilon > 0$ a.s., $\sup_{t_i^n\leq t\leq t_{i+1}^n}|X(t) - X(t_i^n)| \leq (1+\varepsilon)2^n/m$ for all $n$ large and any $i$. The result now follows from easy manipulations. □

Fix $\beta$ and $\varepsilon_1 > 0$ so that $0 < 2(1-\beta) < \varepsilon_1$ and $\varepsilon_1 < \varepsilon_0$ for $\varepsilon_0$ the constant of Corollary 11 (so that Proposition 16 applies to $\varepsilon_1$ and to $1-\beta$).

A difficulty in dealing with a second class particle is in keeping track of its immediate environment. However, in considering how, say, a second class particle at site $x$ at time $t_i^n$ behaves in time interval $[t_i^n, t_{i+1}^n]$, we will be able to deal with the of-order-$2^n$ relevant sites at the same time.

For $x/t_i^n \in \,]1 - 2\lambda + \delta, 1 - 2\rho - \delta[$ and $i = 0, 1, 2, \ldots, m-1$, let $A^{t_i^n}(x)$ be the event that at time $t_i^n$ there is no first class particle occupying site $x$ and that

(23) $$\left|\frac{X^{x,t_i^n}(t_{i+1}^n)}{t_{i+1}^n} - \frac{x}{t_i^n}\right| \geq 2^{-n(1-\beta)}$$

where $(X^{x,t_i^n}(s))_{s\geq t_i^n}$ denotes the position at time $s$ of a (unique) second class particle at site $x$ at time $t_i^n$.

Before analyzing the deviations of these random processes we need a calculus result.

LEMMA 18. *For $1/m < \delta/10$, $\delta < 1/10$, $x \in \,]1 - 2\lambda + \delta, 1 - 2\rho - \delta[$ and $s \in [1 + 1/(2m), 1 + 1/m]$ and*
$$(U_1)'(y) = \begin{cases} \lambda, & \text{if } y \leq 1 - 2\lambda, \\ (1-y)/2, & \text{if } 1 - 2\lambda \leq y \leq 1 - 2\rho, \\ \rho, & \text{if } y \geq 1 - 2\rho, \end{cases}$$

*then the following supremum*
$$U_s(xs) = \sup_{v\in\mathbb{R}}\left\{U_1(v) - (s-1)g\left(\frac{xs-v}{s-1}\right)\right\}$$

is achieved at $v = x$ and for $|v - sx| < 2(s-1)$,
$$U_1(v) - (s-1)g\left(\frac{xs-v}{s-1}\right) \le U_s(xs) - (v-x)^2.$$

PROOF. For any $x$ and $s$ as in the statement, by the definition of function $g$ and the Lipschitz property of $U$,
$$U_s(xs) = \sup_{|v-xs| \le (s-1)} \left\{ U_1(v) - (s-1)g\left(\frac{xs-v}{s-1}\right) \right\}$$
$$= \sup_{|w| < (s-1)} \left\{ U_1(sx) + w\frac{1-xs}{2} - \frac{w^2}{4} - \frac{s-1}{4}\left(1+\frac{w}{s-1}\right)^2 \right\}$$
$$= U_1(xs) + \sup_{|w| \le (s-1)} V(w),$$

where
$$V(w) = \frac{1}{4}\left( (2(1-xs) - w)w - (s-1)\left(1+\frac{w}{s-1}\right)^2 \right),$$
$$V'(w) = \frac{1}{4}\left( 2(1-xs) - 2w - 2\left(1+\frac{w}{s-1}\right) \right)$$
$$= \frac{1}{2}\left( -xs - w\frac{s}{s-1} \right) = -\frac{s}{2}\left( x + \frac{w}{s-1} \right).$$

Thus, as is easily seen $V(w)$ is maximized at $w_0 = -x(s-1)$ or equivalently $U_1(v) - (s-1)g((xs-v)/(s-1))$ is maximized at $v_0 = xs - x(s-1) = x$.

The second derivative of $V$, $V''(w) = -s/(2(s-1))$, which under the condition on $s$ entails that $V''(w) \le -m + 1/2 < -2$. Thus we obtain [on $|w| \le 2(s-1)$]
$$V(w) \le V(w_0) - (w_0 - w)^2;$$
and on $|v - sx| < 2(s-1)$:
$$U_1(v) - (s-1)g\left(\frac{xs-v}{s-1}\right) \le U_1(xs) + V(w_0) - (w_0 - w)^2$$
$$= U_s(xs) - (x-v)^2. \qquad \square$$

PROPOSITION 19. *For $0 < \varepsilon_1 < \varepsilon_0$, event $A^{t_i^n}(x)$ and $t_i^n$ as in (23) and $x/t_i^n \in ]1 - 2\lambda + \delta, 1 - 2\rho - \delta[$,*
$$\mathbb{P}(A^{t_i^n}(x)) \le \exp\left( -\left(\frac{2^n}{m}\right)^{\varepsilon_1/3} \right),$$
*provided that $n$ is sufficiently large.*

PROOF. Consider $x$, as demanded, fixed and assume without loss of generality that the site $x$ is not occupied by a first class particle at time $t_i^n$. We apply Proposition 16 first to the motion of first class particles (so that for this motion at time $t_i^n$ the site $x$ is vacant). Let the associated server process be $z_t(x)$:

$$z_{t_{i+1}^n}\left(\frac{m+i+1}{m+i}(x-2^{n\beta})\right) = \sup_{y\in\mathbb{Z}}\left\{z_{t_i^n}(y) - \xi_{t_{i+1}^n - t_i^n}^{y,t_i^n}\left(\frac{m+i+1}{m+i}(x-2^{n\beta}) - y\right)\right\},$$

where $\xi^{y,t_i^n}(w)$ is $\xi(w)$ derived from the Poisson processes shifted spatially by $y$ and temporally by $t_i^n$. By the definition of supremum we have

$$z_{t_{i+1}^n}\left(\frac{m+i+1}{m+i}(x-2^{n\beta})\right) \geq z_{t_i^n}(x-2^{n\beta}) - \xi_{2^n/m}^{x-2^{n\beta},t_i^n}\left(\frac{x-2^{n\beta}}{m+i}\right).$$

By Proposition 16, using that $2 \cdot 2^{n(1-\varepsilon_1)} \geq 2^{(n+1)(1-\varepsilon_1)} \geq (t_i^n)^{1-\varepsilon_1}$,

$$z_{t_i^n}(x-2^{n\beta}) \geq t_i^n U_1\left(\frac{x-2^{n\beta}}{t_i^n}\right) - 2 \cdot 2^{n(1-\varepsilon_1)}$$

[outside of probability $k\exp(-2^{n\varepsilon_1/2}) \geq k\exp(-(t_i^n)^{\varepsilon_1/2})$], while by Corollary 12 (with $nt = 2^n/m$ and $[nx] = (x-2^{n\beta})/(m+i)$) outside of probability $3\exp(-(2^n/m)^{n\varepsilon_1})$

$$\xi_{2^n/m}^{x-2^{n\beta},t_i^n}\left(\frac{x-2^{n\beta}}{m+i}\right) \leq \frac{2^n}{m}g\left(\frac{x-2^{n\beta}}{2^n/m}\right) + 2 \cdot 2^{n(1-\varepsilon_1)}$$

and so (for $n$ large)

$$z_{t_{i+1}^n}\left(\frac{m+i+1}{m+i}(x-2^{n\beta})\right)$$
$$\geq t_i^n\left(U_1\left(\frac{x-2^{n\beta}}{t_i^n}\right) - \frac{1}{m+i}g\left(\frac{x-2^{n\beta}}{t_i^n}\right)\right) - 4 \cdot 2^{n(1-\varepsilon_1)}$$
$$= t_i^n U_{((m+i+1)/(m+i))}\left(\frac{m+i+1}{m+i}\frac{x-2^{n\beta}}{t_i^n}\right) - 4 \cdot 2^{n(1-\varepsilon_1)},$$

where the last equality holds by Lemma 18 applied with $x = v = (x-2^{n\beta})/t_i^n$ and $s = (m+i+1)/(m+i)$.

Now for $y \geq x$ consider

$$z_{t_i^n}(y) - \xi_{t_{i+1}^n - t_i^n}^{y,t_i^n}\left(\frac{m+i+1}{m+i}(x-2^{n\beta}) - y\right).$$

By Proposition 13 (and Lemma 15), outside probability $4\exp(-(t_i^n)^{\varepsilon_1/2})$

$$z_{t_i^n}(y) \leq t_i^n U_1\left(\frac{y}{t_i^n}\right) + 6 \cdot 2^{n(1-\varepsilon_1)},$$

while by Corollary 12, for all $y$ such that
$$\left| \frac{m+i+1}{m+i}(x - 2^{n\beta}) - y \right| \leq \frac{2^n}{m}$$
we have
$$\xi_{2^n/m}^{y,t_i^n}\left( \frac{m+i+1}{m+i}(x - 2^{n\beta}) - y \right)$$
$$\geq \frac{2^n}{m} g\left( \frac{((m+i+1)/(m+i))(x - 2^{n\beta}) - y}{2^n/m} \right) - 2^{n(1-\varepsilon_1)}$$

outside probability $8 t_i^n \exp(-(2^n/m)^{\varepsilon_1/2})$; thus outside probability $\exp(-2^{2n\varepsilon_1/5})$ for $n$ large,
$$z_{t_i^n}(y) - \xi_{2^n/m}^{y,t_i^n}\left( \frac{m+i+1}{m+i}(x - 2^{n\beta}) - y \right)$$
$$\leq 7 \cdot 2^{n(1-\varepsilon_1)}$$
$$+ t_i^n\left( U_1\left( \frac{y}{t_i^n} \right) - \frac{1}{m+i} g\left( \frac{((m+i+1)/(m+i))(x - 2^{n\beta}) - y}{t_i^n/(m+i)} \right) \right).$$

But by Lemma 18 [with $v = y/t_i^n$, $s - 1 = 1/m + i$ and $x = (x - 2^{n\beta})/t_i^n$ so $v - x = (y - x - 2^{n\beta})/t_i^n \geq 2^{n(\beta-1)}$ since $y \geq x$ by hypothesis; observe also that condition $|v - sx| < 2(s-1)$ is fulfilled] the last term is majorized by
$$t_i^n\left( U_{((m+i+1)/(m+i))}\left( \frac{m+i+1}{m+i} \frac{x - 2^{n\beta}}{t_i^n} \right) - (2^{n(\beta-1)})^2 \right).$$

Thus we obtain, for $n$ large [outside probability $\exp(-2^{2n\varepsilon_1/5}) + \exp(-c2^{n-1}/m)$] that
$$z_{t_{i+1}^n}\left( \frac{m+i+1}{m+i}(x - 2^{n\beta}) \right) > \sup_{y \geq x}\left\{ z_{t_i^n}(y) - \xi_{2^n/m}^{y_1,t_i^n}\left( \frac{m+i+1}{m+i}(x - 2^{n\beta}) - y \right) \right\}.$$

Now consider the exclusion process including the second class particle at site $x$.

Let the height process for this process be denoted $z'_s$ so that
$$z'_{t_i^n}(u) = \begin{cases} z_{t_i^n}(u), & \text{for } u < x, \\ z_{t_i^n}(u) + 1, & \text{for } u \geq x. \end{cases}$$

The above calculations ensure that outside the same probability
$$z'_{t_{i+1}^n}\left( \frac{m+i+1}{m+i}(x - 2^{n\beta}) \right)$$
$$= \sup_{y < x}\left\{ z'_{t_i^n}(y) - \xi_{2^n/m}^{y_1,t_i^n}\left( \frac{m+i+1}{m+i}(x - 2^{n\beta}) - y \right) \right\}$$

$$= \sup_{y<x}\left\{z_{t_i^n}(y) - \xi_{2^n/m}^{y_1,t_i^n}\left(\frac{m+i+1}{m+i}(x-2^{n\beta})-y\right)\right\}$$

$$= z_{t_{i+1}^n}\left(\frac{m+i+1}{m+i}(x-2^{n\beta})\right).$$

This implies that the second class particle must be strictly to the right of site $((m+i+1)/(m+i))(x-2^{n\beta})$ at time $t_{i+1}^n$. In a similar way we obtain outside probability $10\exp(-2^{2n\varepsilon_1/5})$

$$z'_{t_{i+1}^n}\left(\frac{m+i+1}{m+i}(x+2^{n\beta})\right)+1$$

$$= z_{t_{i+1}^n}\left(\frac{m+i+1}{m+i}(x+2^{n\beta})\right)$$

$$= \sup_{y\geq x}\left\{z_{t_i^n}(y) - \xi_{2^n/m}^{y,t_i^n}\left(\frac{m+i+1}{m+i}(x+2^{n\beta})-y\right)\right\}$$

$$> \sup_{y>x}\left\{z_{t_i^n}(y) - \xi_{2^n/m}^{y,t_i^n}\left(\frac{m+i+1}{m+i}(x+2^{n\beta})-y\right)\right\}$$

and that the second class particle must be to the left of $((m+i+1)/(m+i))(x+2^{n\beta})$ at time $t_{i+1}^n$. From this the proposition follows easily. $\square$

PROPOSITION 20. *With probability* 1 *for* $n$ *sufficiently large and for all* $x\in](1-2\lambda+2\delta)2^n,(1-2\rho-2\delta)2^n[$ *which are vacant at time* $2^n$ *for first class particles,*

$$\left|\frac{X^{x,2^n}(t_j^N)}{t_j^N} - \frac{x}{2^n}\right| \leq \sum_{r\geq n} m 2^{-r(1-\beta)} \qquad \forall j \text{ and } N\geq n.$$

PROOF. Let $D^{r,i}$ be the event

$$\left\{\exists y\in](1-2\lambda+\delta)2^r,(1-2\rho-\delta)2^r[,\left|\frac{X^{y,t_i^r}(t_{i+1}^r)}{t_{i+1}^r} - \frac{y}{t_i^r}\right|\geq 2^{-r(1-\beta)}\right\};$$

then we have by Proposition 19 that (for $\varepsilon_1 < \varepsilon_0$ and $r$ sufficiently large)

$$P(D^{r,i}) \leq 2\cdot 2^r \exp\left(-\left(\frac{2^r}{m}\right)^{\varepsilon_1/3}\right)$$

and so

$$P\left(\bigcup_{r\geq n}\bigcup_{i=0}^{m-1} D^{r,i}\right) \leq 2m\sum_{r\geq n} 2^r \exp\left(-\left(\frac{2^r}{m}\right)^{\varepsilon_1/3}\right),$$

which tends to zero as $n$ becomes large. So we have a.s. for $n$ sufficiently large

$$\bigcup_{r \geq n} \bigcup_{i=0}^{m-1} D^{r,i} \qquad \text{does not occur}$$

and also

$$\sum_{r \geq n} m 2^{-r(1-\beta)} \leq \delta.$$

Now notice that if for $r \geq n$ and $i \in \{0, 1, 2, \ldots, m-1\}$,

$$\left| \frac{X^{x,2^n}(t_i^r)}{t_i^r} - \frac{x}{2^n} \right| \leq \sum_{k=n}^{r-1} m 2^{-k(1-\beta)} + i 2^{-r(1-\beta)},$$

then by our choice of $n$, $(X^{x,2^n}(t_i^r))/t_i^r \in\, ]1-2\lambda+\delta, 1-2\rho-\delta[$ and so we may apply Proposition 19 to the process $(X^{x,2^n}(s))_{s \geq t_i^r} = (X^{X^{x,2^n}(t_i^r), t_i^r}(s))_{s \geq t_i^r}$ and conclude that on event $(D^{r,i})^c$,

$$\left| \frac{X^{x,2^n}(t_{i+1}^r)}{t_{i+1}^r} - \frac{x}{2^n} \right| \leq \sum_{k=n}^{r-1} m 2^{-k(1-\beta)} + (i+1) 2^{-r(1-\beta)}.$$

So we have, on the event $\{\bigcup_{r \geq n} \bigcup_{i=0}^{m-1} D^{r,i}\}^c$, by induction on $r$ and $i$, that $\forall r \geq n$ and for each $i \in \{0, 1, 2, \ldots, m-1\}$,

$$\left| \frac{X^{x,2^n}(t_i^r)}{t_i^r} - \frac{x}{2^n} \right| \leq \sum_{k=n}^{r-1} m 2^{-k(1-\beta)} + i 2^{-r(1-\beta)}$$

$$\leq \sum_{k=n}^{\infty} m(2^{-k(1-\beta)}) < \delta. \qquad \square$$

PROOF OF THEOREM 1. From the preceding result we have that a.s. for $n$ large for every $x \in\, ](1-2\lambda+2\delta)2^n, (1-2\rho-2\delta)2^n[$

$$\limsup_{r \to \infty} \frac{X^{x,2^n}(t_j^r)}{t_j^r} - \liminf_{r \to \infty} \frac{X^{x,2^n}(t_j^r)}{t_j^r} \leq 2\delta.$$

From Lemma 17 we conclude that a.s. there is an $n$ so that for every $x \in\, ](1-2\lambda+2\delta)2^n, (1-2\rho-2\delta)2^n[$

$$\limsup_{t \to \infty} \frac{X^{x,2^n}(t)}{t} - \liminf_{t \to \infty} \frac{X^{x,2^n}(t)}{t} \leq \frac{12}{5}\delta < 3\delta.$$

Therefore we conclude that the probability that

$$\limsup_{t \to \infty} \frac{X(t)}{t} - \liminf_{t \to \infty} \frac{X(t)}{t} > 3\delta$$

is bounded by

$$\liminf_{n\to\infty} P(X(2^n) \notin \,](1-2\lambda+2\delta)2^n, (1-2\rho-2\delta)2^n[\,) \le 4\delta/(\lambda-\rho).$$

Since $\delta$ is arbitrary, we conclude that

$$U = \lim_{t\to\infty} \frac{X(t)}{t} \qquad \text{exists a.s.}$$

and, since $X(t)/t$ converges in distribution to $U[0,1]$, we easily obtain $U \sim U[0,1]$. $\square$

**Acknowledgments.** We thank the referees and the editor for helpful comments and constructive criticisms that resulted in an improvement of the presentation.

## REFERENCES

[1] DEMBO, A. and ZEITOUNI, O. (1998). *Large Deviations Techniques and Applications*, 2nd ed. Springer, New York. MR1619036
[2] EVANS, L. C. (1998). *Partial Differential Equations.* Amer. Math. Soc., Providence, RI. MR1625845
[3] FERRARI, P. A. (1992). Shock fluctuations in asymmetric simple exclusion. *Probab. Theory Related Fields* **91** 81–101. MR1142763
[4] FERRARI, P. A. and KIPNIS, C. (1995). Second class particles in the rarefaction fan. *Ann. Inst. H. Poincaré Probab. Statist.* **31** 143–154. MR1340034
[5] FERRARI, P. A., KIPNIS, C. and SAADA, E. (1991). Microscopic structure of travelling waves in the asymmetric simple exclusion process. *Ann. Probab.* **19** 226–244. MR1085334
[6] JOHANSSON, K. (2000). Shape fluctuations and random matrices. *Comm. Math. Phys.* **209** 437–476. MR1737991
[7] KESTEN, H. (1993). On the speed of convergence in first-passage percolation. *Ann. Appl. Probab.* **3** 296–338. MR1221154
[8] KIPNIS, C. and LANDIM, C. (1999). *Scaling Limits of Interacting Particle Systems.* Springer, Berlin. MR1707314
[9] LIGGETT, T. M. (1985). *Interacting Particle Systems.* Springer, New York. MR776231
[10] REZAKHANLOU, F. (1991). Hydrodynamic limit for attractive particle systems on $\mathbb{Z}^d$. *Comm. Math. Phys.* **140** 417–448. MR1130693
[11] REZAKHANLOU, F. (1995). Microscopic structure of shocks in one conservation laws. *Ann. Inst. H. Poincaré Anal. Non Linéaire* **12** 119–153. MR1326665
[12] ROST, H. (1981). Nonequilibrium behaviour of a many particle process: Density profile and local equilibria. *Z. Wahrsch. Verw. Gebiete* **58** 41–53. MR635270
[13] SEPPÄLÄINEN, T. (1998). Hydrodynamic scaling, convex duality and asymptotic shapes of growth models. *Markov Process. Related Fields* **4** 1–26. MR1625007
[14] SEPPÄLÄINEN, T. (1999). Existence of hydrodynamics for the totally asymmetric simple $K$-exclusion process. *Ann. Probab.* **27** 361–415. MR1681094
[15] SEPPÄLÄINEN, T. (2001). Hydrodynamic profiles for the totally asymmetric exclusion process with a slow bond. *J. Statist. Phys.* **102** 69–96. MR1819699

[16] SEPPÄLÄINEN, T. (2001). Second class particles as microscopic characteristics in totally asymmetric nearest-neighbor $K$-exclusion processes. *Trans. Amer. Math. Soc.* **353** 4801–4829. MR1852083

CHAIRE DE PROCESSUS STOCHASTIQUE  
INSTITUT DE MATHÉMATIQUES  
ÉCOLE POLYTECHNIQUE FÉDÉRALE DE LAUSANNE  
1015 LAUSANNE  
SWITZERLAND  
E-MAIL: thomas.mountford@epfl.ch  

LAB. TIMC–TIMB  
ENSIMAG–INP GRENOBLE  
PAVILLON D—FACULTÉ DE MÉDECINE  
38706 LA TRONCHE CEDEX  
FRANCE  
E-MAIL: Herve.Guiol@imag.fr